\documentclass[a4paper,reqno, 11pt]{amsart}
\usepackage{graphicx,color}
\usepackage{footmisc}
\usepackage{amsmath}
\usepackage{tikzsymbols}
\usepackage{amsfonts, dsfont}
\usepackage{amssymb}
\usepackage{bbm}
\usepackage{epstopdf}
\usepackage{color}
\numberwithin{equation}{section}

\let\a=\alpha \let\b=\beta    \let\d=\delta \let\e=\varepsilon
       \let\l=\lambda
    \let\n=\nu         \let\p=\pi    
 \let\t=\tau    
   \let\o=\omega
   \let\L=\Lambda

\newcommand{\ip}{\int_{[-\pi,\pi]^2}\frac{d p}{(2\pi)^2}}
\newcommand{\ik}{\int_{[-\pi,\pi]^2}\frac{d k}{(2\pi)^2}}
\newcommand{\ikp}{\int_{[-\pi,\pi]^2}\frac{d k'}{(2\pi)^2}}

\newtheorem{Theorem}{Theorem}

\newtheorem{Remark}{Remark}

\newtheorem{Lemma}{Lemma}

\def\nn{\nonumber}

\def\\{\hfill\break}
\def\={:=}

\let\0=\noindent

\def\tende#1{\,\vtop{\ialign{##\crcr\rightarrowfill\crcr\noalign{\kern-1pt
    \nointerlineskip} \hskip3.pt${\scriptstyle #1}$\hskip3.pt\crcr}}\,}
\def\otto{\,{\kern-1.truept\leftarrow\kern-5.truept\to\kern-1.truept}\,}

\def\to{\rightarrow}

\def\qed{\hfill\raise1pt\hbox{\vrule height5pt width5pt depth0pt}}

\def\ul#1{{\underline#1}}

\def\be{\begin{equation}}
\def\ee{\end{equation}}
\def\bea{\begin{eqnarray}}
\def\eea{\end{eqnarray}}
\def\nn{\nonumber}

\def\bml{\begin{multline}}
\def\eml{\end{multline}}

\begin{document}
\title{Non-integrable dimer models: universality and scaling relations}
\author{Alessandro Giuliani}
\address{Dipartimento di Matematica e Fisica Universit\`a di Roma Tre \\ {
{L.go S. L. Murialdo 1, 00146 Roma, Italy}}}
\email{giuliani@mat.uniroma3.it}
\author{Fabio Lucio Toninelli}
\address{Univ Lyon, CNRS,  Universit\'e Claude Bernard Lyon 1\\
{UMR 5208, Institut Camille Jordan,
69622 Villeurbanne cedex, France}}
\email{toninelli@math.univ-lyon1.fr}

\begin{abstract} 
  In the last few years, the methods of constructive Fermionic
  Renormalization Group have been successfully  applied to the study of
  the scaling limit of several two-dimensional statistical mechanics
  models at the critical point, including: weakly non-integrable 2D Ising
  models, Ashkin-Teller, 8-Vertex, and close-packed interacting dimer
  models.  In this note, we focus on the illustrative example of
  the interacting dimer model and review some of the universality
  results derived in this context. In particular, we discuss the
  massless Gaussian free field (GFF) behavior of the height
  fluctuations. It turns out that GFF behavior is connected with a
  remarkable identity (`Haldane' or 'Kadanoff relation') between an
  amplitude and an anomalous critical exponent, characterizing the
  large distance behavior of the dimer-dimer correlations.
\end{abstract}

\maketitle

\footnotetext{\copyright\, 2019 by the authors. This paper may be reproduced, in its
entirety, for non-commercial purposes. }

\section{Introduction}

The scaling limit of the Gibbs measure of a statistical mechanics model at a second order phase transition is expected to be universal, in particular, to be robust under `irrelevant' perturbations and  
conformally invariant. Conceptually, the route towards universality is clear: one should first integrate out the small-scale degrees of freedom, then rescale the variables associated with the large-scale fluctuations, and show 
that the critical model reaches a fixed point, under iterations of this procedure (`Wilsonian' Renormalization Group (RG)). On general grounds, the fixed point is expected to be conformally invariant: therefore, in a second step, one can use 
the methods of Conformal Field Theory (CFT) to classify and characterize all the possible conformally invariant fixed points (CFT methods provide, in fact, a complete classification of such fixed 
point theories in two dimensions; remarkably, there has been recent progress in the characterization of three dimensional conformally invariant theories \cite{PRV}, too). The conformal fixed point theory of interest for the description of a given statistical mechanics system can often be identified 
by using specific information on the critical exponents, typically available from the RG construction. 

Even though widely accepted and believed to be correct, there are only
few cases, mostly in two dimensions, for which this procedure and/or
its predictions  can be rigorously confirmed.
%The development of a consistent mathematical theory of universality is one of the central challenges of mathematical physics, of key importance for 
%An key challenge for mathematical physics
%In this sense, developing a consistent mathematical theory of universality is one of the most 
%the foundations of universality theory are still in a primitive stage; 

\begin{enumerate}
\item A first class of critical models for which the existence and
  conformal invariance of the scaling limit is rigorously known
  consists of two-dimensional, integrable Ising and dimer models on
  isoradial graphs \cite{CHI, Duminil-Copin, Dubedat, Kenyon, Smirnov}. The key technical tool used to prove
  conformal invariance is discrete holomorphicity, which is a
  manifestation of integrability \cite{AB}. The method is
  flexible enough to prove robustness of the scaling limit under
  geometric deformations (e.g., of the domain, or of the underlying
  lattice), but it is not able to explain universality under
  perturbations of the microscopic Hamiltonian. [For a proof of
  conformal invariance of crossing probabilities via discrete
  holomorphicity methods in a non-integrable model, see
  \cite{Smirnov_percolation}.]
\item A second class of critical models for which several predictions of Wilsonian RG and CFT have been rigorously substantiated consists of { non-integrable} perturbations of determinantal models, such as
interacting dimers \cite{GMT17a,GMT17b,GMT19}, Ashkin-Teller and vertex models \cite{BFM, F, GM, GMT19}. These results are based on a constructive, fermionic, version of Wilsonian RG: they allow 
to construct the bulk scaling limit of `local fermionic observables' and to prove scaling relations among critical exponents, 
but they are not flexible enough yet to compute `non-local' observables (such as spin-spin correlations in perturbed Ising or monomer-monomer correlations in perturbed dimers) or to 
accommodate geometric deformations of the domain. [For a first result on perturbed Ising models, based on probabilistic tools, see \cite{ADTW}.] 
\end{enumerate}

Not much is rigorously known about the existence and nature of the
scaling limit of other critical models in two and three dimensions. It
is a central challenge of mathematical physics for the coming years
to extend and effectively combine the available techniques, in order
to cover new models of physical interest and new instances of
universality. In this paper, we review a few selected results on the
{universality} of {non-integrable} two dimensional models, based on
the fermionic RG methods mentioned above. For definiteness, we
restrict our attention to `interacting dimer models'. We first
introduce them informally, and give a first overview both of the
`classical' known results, and of our new results.  Next, in the
following sections, we introduce their definitions and state the
relevant results on their critical behavior in a more precise way.

\subsection{Model and  results: an overview.}

\begin{figure}
\centering\includegraphics[width=.8\textwidth]{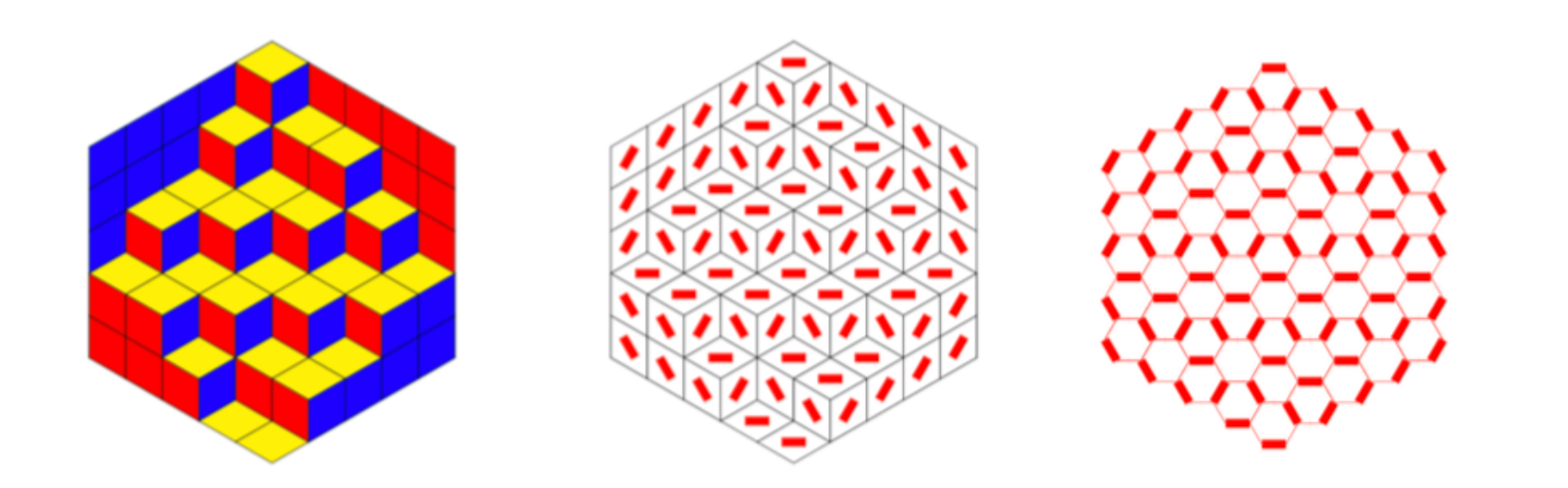}
\caption{A close-packed dimer configuration on a portion of the honeycomb graph (right); the corresponding rhombus tiling picture (center); the tiling can be seen as a stepped monotone surface (left), whose height w.r.t. the horizontal plane is the height function of the dimer configuration.}
\label{fig1}
\end{figure}

Dimer models at close packing on planar, bipartite, graphs are highly
simplified models of {random surfaces} or of liquids of {anisotropic
  molecules} at high density: the connection between these two
apparently unrelated classes of systems is mediated by the notion of
height function, which is in one-to-one correspondence with
close-packed dimer configurations, as illustrated in Fig.\ref{fig1}.

\medskip

{\it Integrable dimer models.} A remarkable feature of close-packed
dimer models is that there is a natural family of exactly solvable,
critical, models, which exhibit a very rich and interesting
behavior. This family is parametrized, on the one hand, by the
specific periodic, infinite, bipartite planar lattice on which the
dimer configuration lives; on the other hand, by the positive weights
(or `activities') associated to the edges of the graph. The exact
solution is based on a determinantal representation of the partition
function, due to Kasteleyn and Temperley-Fisher
\cite{K1,Kenyonnotes,TF}, valid for generic edge weights.  The edge
weights control the average slope of the height: in this sense, the
family of solvable dimer models describes discrete random surfaces
with different average slopes (or `tilt').  The solution shows that
there is an open set of edge weights for which the dimer-dimer
correlations decay algebraically and, correspondingly, the height
fluctuations are asymptotically described by a massless Gaussian Free
Field (GFF) at large distances: this critical phase is called `liquid'
or `rough', depending on whether one focuses on the behavior of the
dimer correlations or of the height profiles. This critical phase
displays a subtle form of universality, in that the variance of the
GFF height fluctuations is asymptotically, at large distances, {\it
  independent} of the dimer weights (in particular, of the
average slope of the height profile)  and of the underlying graph.

\begin{Remark}
  In the case that the underlying lattice is the honeycomb lattice, as in Fig.\ref{fig1}, 
  dimer configurations can be seen as stepped interfaces
  of minimal surface area (given the boundary height), i.e., they can be seen as domain
  walls for the 3D Ising model at zero temperature
  \cite{CerfKenyon}. In fact, the probability measure on dimer
  configurations induced by setting all edge weights equal to $1$ is
  the same as the zero-temperature measure of the 3D Ising model with
  suitable \emph{tilted} Dobrushin-type boundary conditions.  By using
  this connection, one recognizes that the GFF nature of the height
  fluctuations proves the existence of a rough phase in 3D Ising at
  $T=0$ with such boundary conditions. It is an
  open problem to prove that the rough phase persists at small,
  positive, temperatures.
\end{Remark}

\bigskip

{\it Non-integrable dimers: main results (in brief).}  Wilsonian RG
and the bosonization method\footnote{See, e.g., \cite[Sect.1.2]{GMT17b} for a brief introduction to bosonization.} 
suggest that the GFF nature of
the height fluctuations should be robust under non-integrable
perturbations of the dimer model. In order to test this prediction in
a concrete setting, we consider a class of `interacting' dimer models,
including the 6-vertex model (in its dimer representation \cite{Ba1, F}) and {non-integrable} variants thereof. For weak enough
interactions (weak but independent of the system size), we prove that
the height fluctuations still converge to {GFF}, as in the integrable
case. However, in this case, the variance appears to depend on the
interaction and on the dimer weights, see Section \ref{sec4}.
Therefore, the form of universality exhibited by the integrable
dimer model seems to break down as soon as we move out of the exactly
solvable case. Remarkably, a new form of universality emerges in the
interacting case: the (pre-factor of the) variance equals the
anomalous critical exponent of the dimer correlations. This is an
instance of a `Kadanoff' or `Haldane' scaling relation, see
\cite[Sect.1]{GMT19} for a brief introduction to these `weak
universality' relations.

\bigskip

In the next sections, we will describe the models of interest and
state our main results more precisely: in Section \ref{sec2} we
introduce the family of solvable dimer models and briefly review a
selection of known results about its correlation functions; in
Section \ref{sec3} we introduce the class of non-integrable dimer models
we consider, state our main results and give a brief sketch of the
ideas involved in the proof (see \cite{GMT19} for a complete proof);
in Section \ref{sec4} we illustrate the Haldane scaling relation by a
first order computation: we compute both the pre-factor of the
variance and the critical exponent at linear order in the perturbation
strength, and check that the two results agree (the computation shows
how non-trivial and remarkable the result is: already at first order
the validity of the Haldane relation requires very subtle
cancellations). Our explicit computations proves, in particular, that
the critical exponent is non-universal, i.e., it depends both on the interaction strength and
on the dimer weights.

\section{Non-interacting dimers}\label{sec2}

\begin{figure}[h]
\centering
\includegraphics[width=6cm]{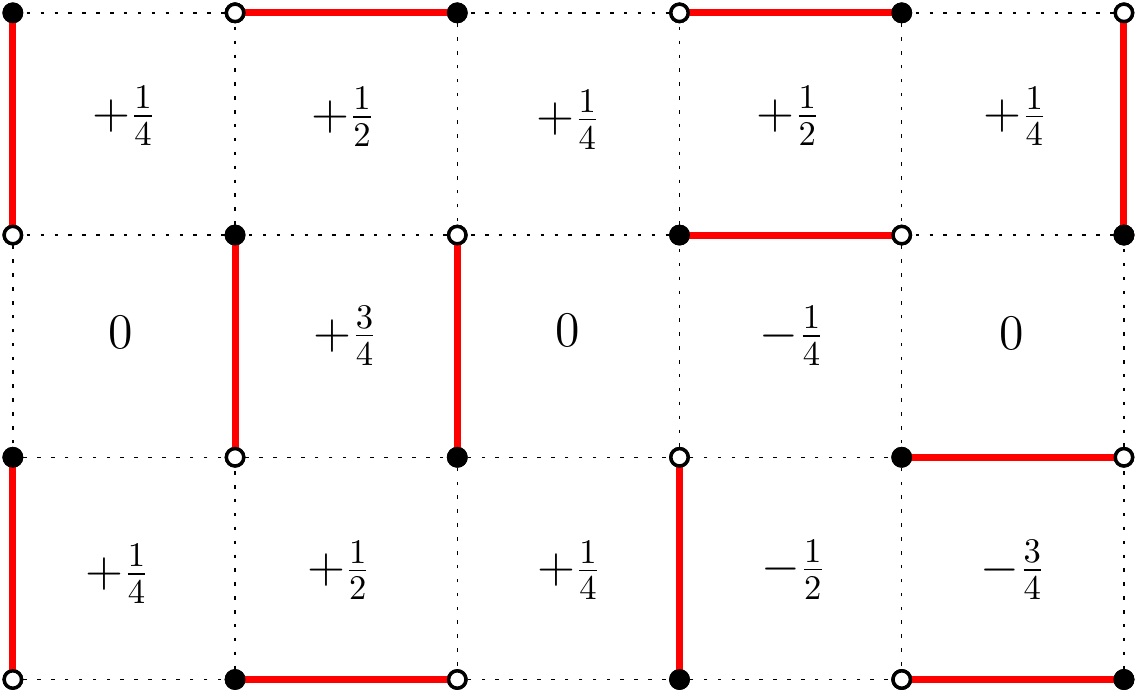}\label{fig2}
\caption{A close-packed dimer configuration on a domain of the square grid and the corresponding height function (the height offset has been fixed to zero at the central face).}
\end{figure}

For simplicity, we restrict to dimers on the square lattice. A
close-packed dimer configuration (or perfect matching) on $\L\subset \mathbb Z^2$ is a
collection of hard rods of length 1, which can be arranged on $\L$ in
such a way that they cover all the vertices of $\L$ exactly once, see
Fig. \ref{fig2}. It is important  that $\L$ is bipartite: we
emphasize this fact by coloring white and black the vertices of the
two sublattices.

Any dimer configuration is in one-to-one correspondence with a height profile, defined on the faces up to an overall additive constant. The height profile is defined by the differences between the values of the height at the face 
$f$ and $f'$:
\begin{equation} \label{hdiff}h(f')-h(f)=\sum_{b\in C_{f\to
      f'}}\sigma_b({\mathds 1}_{b}-1/4),\end{equation} where the sum
runs over the bonds crossed by a lattice path from $f$ to $f'$,
$\sigma_b$ is equal to $+1$ or $-1$, depending on whether $b$ crossed
with the white site on the right or on the left and $\mathds 1_b$ is
the indicator function of having a dimer at $b$. A key point of the
definition of the height function is that the right side is
independent of the choice of the lattice path from $f$ to $f'$
(independence follows from the close-packing condition).

\medskip

The family of integrable, close-packed, dimer models that we informally introduced in the previous section is defined by the following partition function: 
\begin{eqnarray}
  \label{eq:Z0L}
 Z^0_L=\sum_{M\in{\Omega}_L}p_{L,0}(M), \qquad p_{L,0}(M)=\prod_{b\in M}t_{r(b)},  
\end{eqnarray}
defined on a discrete torus $\mathbb T_L$ of side $L$ (see
\cite[Sect.2.1]{GMT19}): here ${\Omega}_L$ is the set of close-packed
dimer configurations on the discrete torus, and $r(b)\in\{1,2,3,4\}$
label the `type' of edge (we let $r=1$ if the edge is horizontal
with the white site on the right, $r=2$ if it is vertical with the
white site on the top, $r=3$ if it is horizontal with the white
site on the left, $r=4$ if it is vertical with the white site on the
bottom).  The model is parametrized by the choice $t_1,t_2,t_3,t_4$;
without loss of generality, we can set $t_4=1$, and we shall do so in
the following. The dimer weights $t_j$ play the role of chemical
potentials, fixing in particular the average slope:
\[\mathbb E_L[h(f+e_i)-h(f)]=\rho_i(t_1,t_2,t_3), \quad i=1,2,\]
where $\mathbb E_L$ indicates the statistical average w.r.t. the  probability measure $\mathbb P_L$ induced by the  weights $p_{L,0}(M)$.

\medskip

As anticipated above, this dimer model is exactly solvable \cite{K1, TF}: for example, 
$$Z^0_L=|\det K|,$$
where $K=K(\ul t)$ is a complex adjacency matrix, known as the {\it Kasteleyn} matrix: its elements are labelled by a pair of sites of different color and are non-zero only if the pair forms a nearest-neighbor edge of the square lattice;
the value of the element of $K$ corresponding to an edge of type $1,2,3,4$ is, correspondingly, 
\begin{equation} K_1=t_1,\quad K_2=it_2,\quad K_3=-t_3,\quad K_4=-i.\label{Kj}\end{equation}

Remarkably, also the multipoint dimer correlation functions can be
explicitly computed. For instance, if $b(x,j)$ is the bond of type
$j$ and black site $x$,  the two-point dimer 
correlation reads (in the special case of two dimers of type 1;
similar formulas are valid for the other cases):
$$\mathbb E[ {\mathds 1}_{b(x,1)};{\mathds 1}_{b(y,1)}
]=-t_1^2\,K^{-1}(x,y)\,K^{-1}(y,x),$$ where: $\mathbb E$ is the
expectation w.r.t. $\mathbb P$, the weak limit of $\mathbb P_L$ as $L\to\infty$; the semi-colon indicates `truncation' (i.e., $\mathbb E(A;B)=\mathbb E(A;B)-\mathbb E(A)\,\mathbb E(B)$); and
\begin{eqnarray}\label{K-1} && K^{-1}(x,y)=\ik\frac{e^{-i k({ x}-{ y})}}{\mu(k)}, \nonumber\\
 {\rm with}  &&  \mu(k)=t_1+ it_2 e^{i k_1}-t_3e^{i k_1+i k_2}-i e^{i k_2}.\nonumber\end{eqnarray}
Note that the zeros of $\mu(k)$ lie at the intersection of two circles in the complex plane, i.e., they are defined by the equation
$$  e^{i k_2}=\frac{t_1+i t_2 e^{i k_1}}{i+ t_3 e^{i k_1}}.$$
`Generically', these two circles intersect transversally in two points: in this situation $K^{-1}(x,y)$ decays to zero algebraically, like $|x-y|^{-1}$, as $|x-y|\to\infty$; this algebraic decay is often referred to by saying that the 
model is {\it critical}.

\medskip

Once the two-point dimer correlations are explicitly known, one can compute the fluctuations of the height difference. In particular, the variance  w.r.t. $\mathbb P$ 
of the height difference diverges logarithmically:
\begin{equation} \label{var}
{\rm Var}_{\mathbb P}[h(f)-h(f')]= \frac{1}{\pi^2}\log |f-f'|+O(1)
\end{equation}
as the distance $|f-f'|$ tends to infinity \cite{Kenyonnotes, KOS}. 

\medskip

\begin{Remark}
The pre-factor $1/\pi^2$ in front of the logarithm is {\it independent} of $t_1,t_2,t_3$: this `universality' property is not accidental and is in fact related to a maximality (`Harnack') property of 
the spectral curve of integrable dimer models \cite{KOS}.
  
\end{Remark}

\medskip

The proof of \eqref{var} is not trivial; it is based on a sufficiently
smart combination of the following ingredients (the `smart' part is in
the use of the third ingredient): (1) the definition of height
difference, see \eqref{hdiff}; (2) the explicitly formula of
$\mathbb E[ \mathds 1_{b};\mathds 1_{b'}]$; (3) the
path-independence of the height (in order to get a sensible expression
in the large distance limit, one first needs to properly deform the
two paths from $f$ to $f'$ involved in the computation of the
variance; next, one can pass to a continuum limit, by replacing
discrete sums with integrals and by replacing the finite distance
dimer-dimer correlation with its large-distance asymptotic behavior,
up to error terms that can be explicitly estimated \cite{KOS}).

\medskip

Building upon these ingredients one can refine \eqref{var} in various directions, in particular one can prove \cite{Kenyon} that: 
\begin{itemize} \item height fluctuations converge to a massless GFF after proper coarse graining and rescaling, 
\item the scaling limit of the height field is conformally covariant.
\end{itemize}

It is very natural to ask whether these features are robust under perturbations of the Gibbs measure that break the determinant structure of Kasteleyn's solution. This point is discussed in the following section. 

\section{Interacting dimers: model and main results}\label{sec3}

We consider a family of `interacting' dimer models, defined by the following partition function: 
\begin{eqnarray}
  \label{eq:ZlL}
  Z^\l_L=\sum_{M\in{\Omega}_L}p_{L,\l}(M), \qquad p_{L,\l}(M)=\Big(\prod_{b\in M}t_{r(b)}\Big)e^{{\l} \sum_{x\in \L}{f}(\t_x M)},
\end{eqnarray}
where: ${\l}$ is the interaction strength (to be thought of as
`small'), $\L$ is the lattice of black sites in $\mathbb T_L$, ${f}$
is a local function of the dimer configuration around the origin, and
$\t_x$ is the `translation operator' by the lattice vector $x$. In
analogy with the non-interacting case, we let $\mathbb P_{L,\l}$ be 
the finite volume Gibbs measure associated with weights $p_{L,\l}$ and
$\mathbb P_\l$ its infinite-volume limit (existence of the limit is part of the results in \cite{GMT19}).

\medskip

\begin{Remark}
For a special choice of $f$, the model reduces to the 6-vertex model, which is integrable by Bethe ansatz (but the solution is not determinantal), see \cite[eq.(2.15)]{GMT19}. However, generically, the model is non-integrable.  
\end{Remark}
 
\begin{figure}
	\begin{center}
		\includegraphics[height=2.4cm]{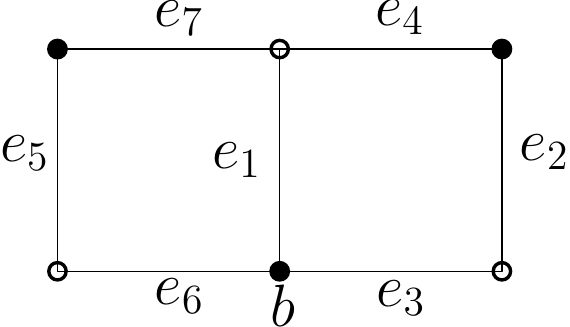}
		\caption{The edges appearing in \eqref{eq:plaq}. $b$ is any fixed black vertex, say the one of coordinates $(0,0)$.}
		\label{fig:esempi}
	\end{center}
\end{figure}

Possibly, the simplest choice of $f$ that makes the model non-integrable is the `plaquette interaction' previously considered in \cite{A,GMT17a,GMT17b}: 
 \begin{equation}
    \label{eq:plaq}
f=\mathds 1_{e_1}\mathds 1_{e_2}+\mathds 1_{e_3}\mathds 1_{e_4}+\mathds 1_{e_1}\mathds 1_{e_5}+\mathds 1_{e_6}\mathds 1_{e_7}
  \end{equation}
where $e_1,\dots,e_7$ are the edges in Fig. \ref{fig:esempi}. 
Our main results do not depend on the specific choice of $f$. They concern the asymptotic behavior of dimer-dimer correlations and of the fluctuations of the height difference and are stated next. 

\medskip

\begin{Theorem}
  \cite{GMT19} \label{th:1} Let $t_1,t_2,t_3$ be such that $\mu(k)$
  has two distinct simple zeros, $p^\pm$ (in particular, the ratio of
  $\a_\o:=\partial_{k_1}\mu(p^\o)$ and
  $\b_\o:=\partial_{k_2}\mu(p^\o)$ is not real). Then, for $\l$ small
  enough, \bea&& \mathbb E_\l[\mathds 1_{b(x,r)};\mathds 1_{b(0,r')}]=
  {\frac1{4\pi^2} \sum_{\o=\pm} \frac{\bar K_{\o, r}\bar K_{\o, r'}}{(\bar \b_\o x_1-\bar \a_\o x_2)^2}}\label{dimdim}\\
  &&{+\frac1{4\pi^2} \sum_{\o=\pm} \frac{\bar H_{-\o, r} \bar H_{\o,
        r'} }{|\bar \b_\o x_1-\bar \a_\o x_2|^{2\n} } e^{-i(\bar
      p^\o-\bar p^{-\o}) x}}+{\bar R_{r,r'}(x)} \;,\nonumber\eea
  where\footnote{A warning on notation: given a quantity (such as
    $\alpha_\o,p^\o$) referring to the non-interacting model, the
    corresponding $\lambda$-dependent quantity for the interacting
    model will be distinguished by a bar, such as $\bar\alpha_\o,$
    etc. On the other hand, we denote by $z^*$ the complex conjugate
    of a number $z$.}: $\bar K_{\o, r}$, $\bar H_{\o,r},$
  $\bar \a_\o$, $\bar \b_\o$, $\bar p^\o$, $\n$ are analytic functions
  of $\l$, such that $\bar \a_\o\big|_{\l=0}=\a_\o$,
  $\bar \b_\o\big|_{\l=0}=\b_\o$, $\bar p^\o\big|_{\l=0}=p^\o$ and
  $\bar K_{\o, r}\big|_{\l=0}=\bar H_{\o,r}\big|_{\l=0}=K_r e^{-i p^\o
    v_r}$. Recall that $K_r$ were defined in \eqref{Kj}; moreover,
  \begin{eqnarray}
    \label{eq:vr}
v_1=(0,0), v_2=(-1,0), v_3=(-1,-1), v_4=(0,-1).    
  \end{eqnarray}
These constants satisfy the following symmetry relations:
\begin{eqnarray}
  \label{eq:symmab}
&&  \bar \alpha_\o^*=-\bar\alpha_{-\o},\hskip.75truecm 
  \bar\beta_\o^*=-\bar\beta_{-\o} \\
  \label{eq:symmK}
 && \bar K_{\o, r}^*=\bar K_{-\o,r}  ,\quad 
  \bar H_{\o, r}^*=\bar H_{-\o,r},\\ 
  \label{eq:symmp}
  &&\bar p^++\bar p^-=(\pi,\pi).
\end{eqnarray}

  Moreover, $$\nu(\l)=1+a\l+O(\l^2)$$ and, generically, $a\neq
  0$. Finally, $\bar R_{r,r'}(x,x')=O(|x-x'|^{-5/2})$ (the exponent
  $5/2$ could be replaced by any $\d<3$ provided $\l$ is small
  enough).
\end{Theorem}
 
 \bigskip

 A few comments are in order:
 \begin{itemize}
 \item  The proof provides a constructive algorithm for computing $\bar K_{\o,r}, \bar H_{\o,r}$, etc, at any desired precision. However, we do not have closed formulas for these quantities,
and we do not expect that it is possible to obtain any by other methods. 

\item 
 Generically, the {\it anomalous exponent} $\nu$ has a non-zero first order contribution in $\l$: therefore, it is larger or smaller than $1$, depending on whether $a \l$ is positive or negative. 
In particular, the asymptotic, large-distance, behavior of the dimer-dimer correlation is dominated by the first or second term in \eqref{dimdim}, depending on whether $a\l$ is positive or negative. 

\item  Once we have such a refined asymptotics as \eqref{dimdim}, we can plug it into the formula for the height variance, 
$$\mathbb E_\l[{h(f)-h(f');h(f)-h(f')}]=\sum_{b,b'\in C_{f\to f'}}\sigma_b\sigma_{b'}\mathbb E_\l[
\mathds 1_{b};\mathds 1_{b'} ],$$ in order to understand its growth as
$|f-f'|\to \infty$. As anticipated above, one expects that the growth
is still logarithmic, like in the non-interacting case.  However, this
is not completely straightforward: as remarked in the previous item,
if $a\l<0$ the asymptotic behavior of the dimer-dimer correlation is
dominated by the second term in \eqref{dimdim}, which is characterized
by the critical exponent $\nu<1$. Therefore, after double-summation
over $b,b'$, one may fear that the interacting variance
$\mathbb E_\l[{h(f)-h(f');h(f)-h(f')}]$ grows like $|f-f'|^{2(1-\nu)}$ at
large distances, rather than logarithmically.  Remarkably, this is not
the case, thanks to the oscillating factor
$e^{-i(\bar p^\o-\bar p^{-\o}) x}$ in front of the second term of
\eqref{dimdim}.

 \end{itemize}

 \begin{Theorem}  \cite{GMT19}
   \label{th:2}
Under the same hypothesis of Theorem 1, 
\[
{\rm Var}_{\mathbb P_\l}(h(f)-h(f'))= \frac{{A(\lambda)}}{\pi^2}\log |f-f'|+O(1)
\]
where 
$$A(\l)=-\Biggl[\frac{\bar K_{\o,1}+\bar K_{\o,2}}{\bar \b_\o}\Biggr]^2=-\Biggl[\frac{\bar K_{\o,1}+\bar K_{\o,4}}{\bar \a_\o}\Biggr]^2.$$
In general, $A(\l)$ depends explicitly on $\l,t_1,t_2,t_3$ and on the type of interaction potential in \eqref{eq:ZlL}. Moreover, 
\begin{equation} \label{Ha}{A(\l)}={\n(\l)}.\end{equation}
 \end{Theorem}

 The proof of Theorem \ref{th:2} provides two different algorithms for
 computing the coefficients of the analytic functions $A$ and
 $\nu$. At first order in $\l$, one can check directly the validity of
 \eqref{Ha}. The computation is very instructive and, as a by-product,
 it shows that $A(\l)$ depends explicitly on $\l,t_1,t_2,t_3$; see
 Sect.\ref{sec4} for a proof.  However, a direct proof of the validity
 of $A=\nu$ by direct inspection and comparison of the two power
 series defining $A$ and $\nu$ seems hopeless (the computation in the
 next section is a clear indication).  The proof of \eqref{Ha} is
 based on a subtle cancellation mechanism, which uses a comparison
 between the lattice Ward Identities of the dimer model and the
 emergent chiral Ward Identities of a `reference' continuum model (the
 `infrared fixed point theory'), see \cite[Sect.4 and 5]{GMT19}.

\medskip

The identity $A(\l)=\nu(\l)$ is the analogue of one of the scaling relations (the one between compressibility and the density-density critical exponent) proposed by D. Haldane in the context of Luttinger 
liquids \cite{Ha}. 
It is also related to one of the scaling relations among the critical exponents of the 8-vertex and of the Ashkin-Teller model, proposed by L. Kadanoff \cite{Ka}, see the discussion right before eq.(1.2) of \cite{GMT19}. For this reason, we call it `Haldane' or `Kadanoff' relation.

\medskip

 There are not many previous examples of rigorously established Haldane or Kadanoff relations: the known cases are mostly restricted to exactly solved models of interacting fermions or quantum spins (Luttinger model, or the XXZ chain \cite{Ha}) and non-integrable perturbations thereof \cite{BMvarii}. 

\medskip

As stated, Theorem 2 concerns only the asymptotic growth of the variance of the height difference. However, a slight extension of its proof, along the same lines as \cite[Sect.7.3]{GMT17a}, proves that, after 
coarse-graining, $h(f)$ converges to $\phi$, the massless GFF of
covariance 
$$\mathbb E_\l(\phi(x) \phi(y))=-\frac{A(\lambda)}{2\pi^2}\log|x-y|, $$
in the same sense as \cite[Theorem 3]{GMT17a}.

Theorem \ref{th:2} and its extension mentioned in the previous item prove in a very strong and sharp sense that the random surface associated with the interacting dimer model is in a rough phase, characterized 
by logarithmic fluctuations. In this sense, our theorem is of interest in the context of the fluctuation theory of discrete random surface models, which is a classical topic in probability and statistical mechanics. Previous, related, 
results include the proofs of logarithmic height fluctuations in anharmonic crystals,  SOS model, 6-vertex and Ginzburg-Landau type models \cite{BLL, CoSpe, F, FroSpe, Giacomin-Olla-Spohn,  Ioffe-Shlosman-Velenik,Mill, Milos-Peled,  Naddaf-Spencer}.

\bigskip

The proofs of Theorems \ref{th:1} and \ref{th:2} are hard and
lengthy, and we refer the reader to \cite{GMT19} for it. Here we limit
ourselves to mention the main steps and ideas involved in the proof:
\begin{enumerate}
\item The starting point is a representation of the non-interacting
  model, characterized by Kasteleyn's determinantal solution, in terms
  of a free fermionic theory in dimension $d=1+1$.
\item Next, we provide an exact representation of the interacting dimer model in terms of an interacting fermionic theory in $d=1+1$, the interaction being (at dominant order) a quartic fermionic 
interaction. In this sense, the interacting dimer model maps into a sort of `fermionic $\phi^4_2$ theory'. 
\item The fermionic model, into which the original dimer model has
  been mapped, can be analyzed via standard fermionic multiscale
  cluster expansion methods (fermionic constructive RG) due, among
  others, to Gawedzki-Kupiainen \cite{GK}, Lesniewski \cite{Le},
  Benfatto-Gallavotti \cite{BG}, Feldman-Magnen-Rivasseau-Trubowitz \cite{FMRT}.
\item The fermionic RG scheme turns out to be convergent if and only if the %fermionic RG scheme allows us to reduce the control of the large distance behavior of the 
infrared flow of the `running coupling constants', controlled by the so-called beta-function equation, is convergent in the limit of large scales. In order to control the RG flow, we 
compare it with the one of a reference continuum model: some of the ingredients involved in the comparison are the Ward Identities (both for the lattice dimer model and for the continuum reference one), 
the Schwinger-Dyson equation and the non-renormalization of the anomalies for the reference model \cite{BMvarii}. 
\item In order to obtain the fine asymptotics of the dimer-dimer
  correlations, as well as the Haldane relation connecting $A$ and
  $\nu$, we need to compare the asymptotic, emergent, chiral Ward
  Identities of the reference model with the exact lattice Ward
  Identities of the dimer model, following from the local conservation
  law of the dimer number, $\sum_{b\sim x}\mathds 1_b=1$, with the sum
  running over edges incident to $x$. The comparison guarantees that
  the ratio $A/\nu$ is `protected' by symmetry (there is no `dressing'
  or `renormalization' due to the interaction).
\end{enumerate}

\bigskip

We conclude this section by mentioning a few open problems and
perspectives:
\begin{itemize}
\item It would be interesting to study the critical theory in finite domains and, in perspective, prove its conformal covariance. From a technical point of view, going in this direction requires a 
non-trivial extension of the RG multiscale methods to the non-translationally-invariant setting, and a sharp control of the RG flow of the marginal boundary running coupling constants. 
Promising results in this direction have been recently obtained in the context of non-planar critical Ising models in the half-plane \cite{AAG}. 
\item So far, we can only compute the scaling limit of the height
  function after coarse graining, cf. for instance
  \cite[Th. 3]{GMT17a}. It would be very interesting to prove a
  central limit theorem on a more local scale; i.e., computing the
  average of $e^{i \alpha (h(f)-h(f'))}$, instead of computing the
  characteristic function for the height function integrated against a
  test function. Similarly, it would be interesting to compute the
  scaling limit of the monomer-monomer correlations. This problem is
  not merely technical: it is strictly connected with the computation
  of the scaling limit of the spin-spin correlations in non-integrable,
  two-dimensional Ising models, which is currently out of reach of
  the current techniques.
\item While we expect the analogs of Theorems \ref{th:1} and
  \ref{th:2} to be true also for the interacting dimer model on the
  honeycomb lattice of Fig. \ref{fig1}, it is not obvious that the
  same qualitative behavior holds for interacting dimers on general
  $\mathbb Z^2$-periodic bipartite planar graphs, with an elementary
  cell containing more than two vertices. In the non-interacting case,
  using $\mathbb Z^2$ or any other $\mathbb Z^2$-periodic bipartite
  planar graph makes essentially no difference \cite{KOS}. However,
  the effective fermionic theory is different for different periodic
  bipartite lattices: the larger the elementary cell, the larger the
  number of `colors' of the fermionic field associated with the
  fermionic description of the system. In the presence of
  interactions, the number of colors is known to affect the
  qualitative behavior of the system (as an example, compare the behavior of the Luttinger model \cite{ML} with that 
  of the 1D Hubbard model \cite{LW}): depending on it and on the
  specific form of the interaction, the system could display either an
  anomalous Fermi liquid behavior, or it could open a gap, entering
  more exotic quantum phases.  It would be very interesting to see
  whether any of these exotic behaviors can arise in interacting dimer
  models on decorated periodic lattices.
\item Finally, as suggested by the remarks in the introduction, it
  would be nice to see whether the current methods can be used to
  prove the existence of a rough, logarithmically correlated,
  low-temperature phase of the interface of 3D Ising with tilted
  Dobrushin boundary conditions. As observed above, the fluctuations
  of the surface at zero temperature can be mapped exactly into a
  problem of non-interacting dimers on the hexagonal lattice. It
  remains to be seen whether low-temperature fluctuations can be
  effectively described in terms of weakly interacting dimers on the
  hexagonal lattice.
\end{itemize}

\section{First order computation}\label{sec4}

The goal of this section is to verify, at first order in $\l$, the validity of the `Haldane', or `Kadanoff', relation, \eqref{Ha}, and to compute explicitly 
(to fix ideas, for the case of the plaquette interaction \eqref{eq:plaq}) the
constant $a$ in the expansion
\begin{eqnarray}
  \label{eq:expansione}
\nu(\l)=A(\l)=1+a\l+O(\l^2).  
\end{eqnarray}
The result is: \be \label{eqa}a=-\frac{4(t_1 t_3+t_2)}{\pi}\frac{
  \cos(p^+_1)\cos(p^+_2)}{\Im(\alpha_+\beta_-)},\ee where we recall
that $p^\pm$ are the two zeros of $\mu(k)$, which are assumed to be
non-degenerate, $\a_\pm=\partial_{k_1}\mu(p^\pm)$ and
$\b_\pm=\partial_{k_2}\mu(p^\pm)$. It is straightforward to check that
the right side of \eqref{eqa} explicitly depends on the dimer weights:
e.g., in the simple case $t_1=t_3=t$, $t_2=1$, in which $p^+=0$,
$p^-=(\pi,\pi)$, one has $a=-2(1+t^2)/(\pi t)$; in the case
$t_1=t_2=t_3=1$, the result coincides with the one found in
\cite[App. B]{GMT17b}, $a=-4/\pi$.

\medskip

We will use the same notations and conventions as in
\cite{GMT19,GMT17a}; we do not repeat here the definitions and we
assume the reader has familiarity in particular with \cite[Sec. 2 and
3]{GMT19} and \cite[App. B]{GMT17b}.  We recall just that
$\mathcal E_0$ denotes the Grassmann Gaussian measure with propagator (cf. \eqref{K-1})
\begin{eqnarray}
  \label{eq:E0}
  g(x-y):= \mathcal E_0(\psi^-_x\psi^+_y)=K^{-1}(x,y)
\end{eqnarray}
For the
following, it is convenient to introduce the rescaled coupling
constant
\begin{eqnarray}
  \label{eq:u}
  u:=(t_1 t_3+t_2)\l.
\end{eqnarray}

\subsection{First-order computation of $A$}

The two-point dimer
correlation function is given as
\begin{eqnarray}
  \label{eq:dd}
  G^{(0,2)}_{r,r'}(x,0):=\mathbb E_\l(\mathds 1_e;\mathds 1_{e'})=\lim_{L\to\infty}\partial_{A_e} \partial_{A_{e'}} \mathcal W_\L(A)|_{A\equiv 0}, 
\end{eqnarray}
where $e$ (resp. $e'$) is the edge of type $r$ (resp. $r'$) with black endpoint of coordinates $x$ (resp. $0$) and $\mathcal W_\L(A)$ is the moment generating function
\begin{eqnarray}
  \label{eq:momgen}
  e^{\mathcal W_\L(A)}:=\sum_{M\in\Omega_L}p_{L,\l}(M)\prod_e e^{A_e \mathds 1_e}.
\end{eqnarray}
Our convention for the Fourier transform of the two-point dimer correlation  is that
\begin{eqnarray}
  \label{eq:Gfourier}
  \hat G^{(0,2)}_{r,r'}(p)=\sum_x e^{-i p x}G^{(0,2)}_{r,r'}(x,0).
\end{eqnarray}
Since we are interested in the long-distance behavior of correlations, we will look at the small-$p$ behavior and in particular at the discontinuity of $\hat G^{(0,2)}_{r,r'}(p)$ at $p=0$.
We recall that, since we are working on the torus, $\exp(\mathcal W_\L(A))$ can be written as the linear
combination of four Grassmann integrals with non-quadratic action,
cf. \cite[Eq. (3.16)]{GMT19}. For the computation of correlations in
the $L\to\infty$ limit and at finite order in perturbation theory, we
can safely replace \cite[Eq. (3.16)]{GMT19} with an expression
involving \emph{a single} Grassmann integral:
\begin{eqnarray} 
  \label{eq:piusemp}
   e^{\mathcal W_\L(A)}=\int D\psi e^{S(\psi)+V(\psi,A)},
\end{eqnarray}
with $S$ and $V$ as in \cite[Eq. (3.16) and (3.17)]{GMT19}.
In the case of the plaquette
interaction \eqref{eq:plaq}, the potential $V(\psi,A)$
 equals (neglecting terms of order $\l^2$
and higher)
\begin{eqnarray}
  \label{eq:Vprimod}
V(\psi,A)=  -\sum_e(e^{A_e}-1)E_e+\l\sum_{\gamma=\{e_1,e_2\}\subset \Lambda}E_{e_1}E_{e_2}e^{A_{e_1}+A_{e_2}}
\end{eqnarray}
where the second sum runs over pairs of parallel edges $\{e_1,e_2\}$
on the boundary of the same face.  In particular, setting $A\equiv 0$
and using the definition \cite[Eq. (3.13)]{GMT19} for $E_e$ in terms
of the Grassmann variables $\psi^\pm$ associated to the endpoints of
the edge $e$, one finds that the potential is exactly quartic in the
Grassmann fields:
\begin{align}\label{eq:V4}
 V_4(\psi):=  V(\psi,0)&=-u\sum_x\psi^+_x\psi^-_x\left[\psi^+_{x+(0,1)}\psi^-_{x-(1,0)}+\psi^+_{x+(1,0)}\psi^-_{x-(0,1)}\right].
 \end{align}
 From  \eqref{eq:dd}, \eqref{eq:momgen} and \eqref{eq:piusemp} we see
 that the two-point dimer correlation function equals, at first order in
 $\lambda$ and in the $L\to\infty$ limit,
\begin{align}
\label{i1}\nonumber
    G^{(0,2)}_{r,r'}(x,0)&=\mathcal E_0( E_{e};E_{e'})\\&-
  \lambda  [\mathcal E_0( E_{e};I^{(1)}_{0,r'})+     \mathcal E_0(I^{(1)}_{x,r};E_{e'})]+\mathcal E_0(E_{e};E_{e'};V_4)
\end{align}
where $\mathcal E_0(\dots;\dots)$ denotes truncated expectation  and 
   \begin{multline}
     I^{(1)}_{x,r}=\left\{
       \begin{array}{lll}
         K_1 K_3\psi^+_{x}\psi^-_x(\psi^+_{x+(0,1)}\psi^-_{x-(1,0)}+\psi^+_{x+(1,0)}\psi^-_{x-(0,1)}) & \text{ if} & r=1  \\
         K_2 K_4\psi^+_{x}\psi^-_{x+v_2}(\psi^+_{x+(0,1)}\psi^-_{x}+\psi^+_{x-(1,0)}\psi^-_{x-(1,1)})& \text{ if} & r=2  \\
         K_1 K_3\psi^+_{x}\psi^-_{x+v_3}(\psi^+_{x-(1,0)}\psi^-_{x-(1,0)}+\psi^+_{x-(0,1)}\psi^-_{x-(0,1)}) &\text{ if} & r=3  \\
         K_2 K_4\psi^+_{x}\psi^-_{x+v_4}(\psi^+_{x+(1,0)}\psi^-_{x}+\psi^+_{x-(0,1)}\psi^-_{x-(1,1)})&\text{ if} & r=4.
       \end{array}
\right.   \end{multline}

% In this appendix we let $K_1=Ae^{B_2}$, $K_2=i$, $K_3=-Ae^{B_1}$, $K_4=-ie^{B_1+B_2}$, even though the specific form of these $K_r$ is not 
% really important. 
% We work in Fourier space, setting (PROBABILMENTE QUESTA FORMULA STA GIA' PRIMA NEL TESTO)

% \be \label{eq:G02F}
%  \hat G^{(0,2)}_{r,r'}(p)=\sum_x e^{-i p x}G^{(0,2)}_{r,r'}(x,0).
% %  G^{(0,2)}_{r,r'}(x,0)=\int\limits_{[-\p,\p]^2}\frac{dp}{(2\p)^2}e^{i p x}\hat G^{(0,2)}_{r,r'}(p).
% \ee
For $\lambda=0$ one finds from \eqref{i1}, \eqref{eq:E0} and from  Lemma \ref{lemma:formulozzo} below
\begin{align}
\label{eq:mardepanza}
  \left. \hat G^{(0,2)}_{r,r'}(p)\right|_{\lambda=0}&=-K_r K_{r'}\ik \frac{e^{-i k v_{r}-i (k+p)v_{r'}}}{\mu(k)\mu(k+p)}\\\nonumber&=
-\frac i{2\pi}\frac{K_r K_{r'}}
   {\alpha_+\beta_--\alpha_-\beta_+}\sum_{\omega=\pm}\frac{D_{-\omega}(p)}{D_\omega(p)}e^{-i p^\omega(v_r+v_{r'})}+R(p).
\end{align}
Here and in the following, $R(p)$ denotes a function that is
continuous at $p=0$ (the precise value of $R(p)$ can change from line
to line).  In \eqref{eq:mardepanza}, $v_r\in \mathbb Z^2$, $r=1,\dots,4$,
is as in \eqref{eq:vr} while
 \begin{eqnarray}
   \label{eq:Domega}
D_\omega(p)=\a_\o p_1+\b_\o p_2,\quad \o=\pm.
 \end{eqnarray}
Next we compute the first-order contribution
   \begin{equation}
-
  \lambda  [\mathcal E_0( E_{e};I^{(1)}_{0,r'})+     \mathcal E_0(I^{(1)}_{x,r};E_{e'})]
   \end{equation}
   in \eqref{i1}.  As explained at the beginning of
   \cite[Sec. B.1]{GMT17b}, by the fermionic Wick theorem, we can
   replace $I^{(1)}_{x,r}$ by its ``linearization''
   $\bar I^{(1)}_{x,r}$ obtained by contracting in all possible ways
   two of its four $\psi$ fields. For instance, with $g(\cdot)$ as in
   \eqref{eq:E0},
   \begin{multline}
     \label{eq:Ibar1}
     \overline I^{(1)}_{x,1}=K_1K_3\left[-g(v_1)\left(\psi^+_{x+(0,1)}\psi^-_{x-(1,0)}+\psi^+_{x+(1,0)}\psi^-_{x-(0,1)}
     \right)\right.\\+g(v_2)\left(\psi^+_{x+(0,1)}\psi^-_x+\psi^+_x\psi^-_{x-(0,1)}
   \right)\\ \left.+g(v_4)\left(\psi^+_{x+(1,0)}\psi^-_x+\psi^+_x\psi^-_{x-(1,0)}
     \right)
   -2g(v_3)\psi^+_x\psi^-_x \right].
   \end{multline}
   In Fourier space (with the conventions of \cite[Eq. (6.1)]{GMT19} for the
   Grassmann fields in momentum space) one has
 \begin{equation}
     \label{eq:Ibar1F}
     \overline I^{(1)}_{x,r}=\ik\ip e^{i p x}\hat\psi^+_{k+p}W_r(k,p)\hat\psi^-_k
\end{equation}
where
\begin{multline}
\label{eq:W1}
  W_1(k,p)=K_1K_3\left[-g(v_1)\left(e^{i(k_1+k_2+p_1)}+e^{i(k_1+k_2+p_2)}
       \right)\right.\\+g(v_2)\left(e^{i(k_2+p_2)}+e^{i k_2}
   \right) \left.+g(v_4)\left(e^{i(k_1+p_1)}+e^{i k_1}
     \right)
   -2g(v_3) \right].
\end{multline}
Similar formulas hold for $r=2,3,4$, with $W_1(k,p)$ replaced by
$W_r(k,p)$.  One easily checks (for $r=1$ this can be verified
immediately from \eqref{eq:W1}) that
\begin{multline}
  \label{eq:W0}
  W_r(k,0)=2 t_r t_{r+2}W(k),\\
  W(k)=\left[g(v_1)e^{i(k_1+k_2)}-g(v_2)e^{ik_2}-g(v_4)e^{i k_1}
    +g(v_3) \right]\\=
  \ikp\frac{(e^{ik_1}-e^{i k_1'})(e^{ik_2}-e^{i k_2'})}{\mu(k')}
\end{multline}
with the convention that $t_4=1$ and that $t_r:=t_{r\!\!\mod \!4}$ if $r>4$.  Note that
\begin{eqnarray}
  \label{eq:nota0}
  {W(p^+)}^*=W(p^-),
\end{eqnarray}
because $p^++p^-=(\pi,\pi)$ and $g(v_1),g(v_2)\in \mathbb R$, while
$g(v_2),g(v_4)\in i\mathbb R$.  Using \eqref{eq:W0} and
\eqref{eq:I0F}, together with \begin{eqnarray}
     \label{eq:I0F}
     E_{e}=K_r\ik \ip e^{i p x}\hat \psi^+_{k+p}\hat\psi^-_ke^{-i k v_r},
                              \end{eqnarray}
                              (if $e$ is, as above, of type $r$ and
with black vertex of coordinates $x$)
we see that
\begin{multline}
\label{eq:vestooss}
-
  \lambda  [\mathcal E_0( E_{e};I^{(1)}_{0,r'})+     \mathcal E_0(I^{(1)}_{x,r};E_{e'})](p)
%  \ik \frac1{\mu(k)\mu(k+p)}[K_r W_{r'}(k+p,-p)e^{-i k v_r}\\+K_{r'}W_r(k,p)e^{-i (k+p)v_{r'}}]
\\=
2 \l\ik \frac1{\mu(k)\mu(k+p)}\times\\
\times[K_r t_{r'}t_{r'+2}e^{-i k v_r}+K_{r'}t_r t_{r+2}e^{-i kv_{r'}}]W(k)
+R(p).
\end{multline}
Thanks to Lemma \ref{lemma:formulozzo}, we can rewrite \eqref{eq:vestooss} as 
\begin{multline}
\label{eq:vestooss2}
\frac{i\l}{\pi(\alpha_+\beta_--\alpha_-\beta_+)}\times
\\
\nonumber\times
\sum_{\o=\pm}\frac{D_{-\o}(p)}{D_\o(p)}\left[
K_r t_{r'}t_{r'+2}e^{-i p^\o v_r}+K_{r'} t_r t_{r+2}e^{-i p^\o v_{r'}}
\right]W(p^\o)
+R(p).
\end{multline}

It remains to compute the term $\mathcal E_0(E_{e};E_{e'};V_4)$ in
\eqref{i1}.  Applying Wick's theorem, we see that either $N=0$ or
$N=2$ of the four fields of $V_4$ are contracted among themselves, the
remaining $4-N$ ones being contracted with fields from $E_{e}$ or
$E_{e'}$.  Here we compute the contribution, call it
$a^{(0,2)}_{r,r'}(x,0)$, from those terms where $N=2$. The
contribution $b^{(0,2)}_{r,r'}(x,0)$ from the terms with $N=0$ is
computed later.  We note that
\begin{eqnarray}
\label{N2}
  \mathcal E_0( E_{e};E_{e'})+\lambda a^{(0,2)}_{r,r'}(x,0)=
   \mathcal E_\l( E_{e};E_{e'})+O(\lambda^2)
\end{eqnarray}
where $\mathcal E_\l(\cdot)$ is the Gaussian Grassmann measure
where the action $S(\psi)=-\sum_e E_e$ has been replaced by
$S(\psi)+ \overline V_4$, with $\overline V_4$ the
linearization of $V_4$ (i.e. the bilinear operator obtained by
contracting in every possible way two of the four fields of $V_4$, as above).  In
our case,
\begin{gather}
  \bar V_4=-2u\sum_z[-g(v_1)\psi^+_{z}\psi^-_{z-(1,1)}+g(v_2)\psi^+_{z}\psi^-_{z-(0,1)}
\\-g(v_3)\psi^+_z\psi^-_z+g(v_4)\psi^+_z\psi^-_{z-(1,0)}].
\end{gather}
Then, we see that the measure $\mathcal E_\l( \cdot)$ is nothing but
the Gaussian Grassmann measure where $\mu(k)$ is replaced by
\begin{equation}
  \bar \mu(k):=\mu(k)-2 u W(k)
\end{equation}
and $W(k)$ was defined in \eqref{eq:W0}.  As already mentioned in
Theorem \ref{th:1}, the ratio of complex numbers $\alpha_\o,\b_\o$ is
not real \cite{GMT19}; therefore we can write uniquely $W(p^\omega)$
as
\begin{eqnarray}
  \label{eq:fpo}
  W(p^\omega)=c^\omega_1\alpha_\omega+c^\omega_2\beta_\omega \text{ with } c^\omega_1,c^\omega_1\in \mathbb R.
\end{eqnarray}
Via Taylor expansion, we have then (always at first order in $\l$)
\begin{eqnarray}
  \label{eq:Taylor}
  \bar \mu(k)=\bar\alpha_\omega(k_1-\bar p^\omega_1)+\bar\beta_\omega(k_2-\bar p^\omega_2)+O(|k-\bar p^\omega|^2)
\end{eqnarray}
where the zeros of $\bar\mu(\cdot)$ are
\begin{eqnarray}
  \label{eq:ptilde}
  \bar p^\omega_j=p^\omega_j+2u c^\omega_j, \quad j=1,2
\end{eqnarray}
and
\begin{multline}
  \label{a+b+}
  \bar  \alpha_\omega=\alpha_\omega+2u[ -\partial_{k_1}W(p^\omega)+c^\omega_1\partial^2_{k_1}\mu(p^\omega)+c^\omega_2\partial^2_{k_1 k_2}\mu(p^\omega)]\\
  \bar  \beta_\omega=\beta_\omega+2u[- \partial_{k_2}W(p^\omega)+c^\omega_2\partial^2_{k_2}\mu(p^\omega)+c^\omega_1\partial^2_{k_1 k_2}\mu(p^\omega)].
\end{multline}
Explicitly,
\begin{align}
  \nn
  \partial_{k_1}W(p^\omega)&=  i [g(v_1)e^{i (p^\o_1+p^\o_2)}-g(v_4)e^{i p^\o_1}]\\
  \nn
  \partial_{k_2}W(p^\omega)&=  i [g(v_1)e^{i (p^\o_1+p^\o_2)}-g(v_2)e^{i p^\o_2}]\\
  \label{adm}
  \partial^2_{k_1}\mu(p^\omega)&=i\alpha_\o, \quad  \partial^2_{ k_2}\mu(p^\omega)=i\beta_\o, \quad  \partial^2_{k_1k_2}\mu(p^\omega)=t_3 e^{i(p^\o_1+p^\o_2)}.
 % e^{B_1+B_2+i p^\o_2}+A
%  e^{B_1+i(p^\o_1+p^\o_2)}.
\end{align}
This implies in particular the symmetry \eqref{eq:symmab} at first
order in $\l$.  Also, from \eqref{eq:nota0} it follows that
$c^\o_j=-c^{-\o}_j$ so that \eqref{eq:symmp} holds a first order.

Altogether, \eqref{N2} equals in Fourier space, for small $p$ and
disregarding $O(\l^2)$ terms,
\begin{eqnarray}
   \label{eq:ordine0dress}
-\frac i{2\pi}\frac{K_r K_{r'}}
   {\bar \alpha_+\bar\beta_--\bar \alpha_-\bar \beta_+}\sum_{\o=\pm}\frac{\bar D_{-\omega}(p)}{\bar D_\omega(p)}e^{-i \bar p^\omega(v_r+v_{r'})}+R(p)
\end{eqnarray}
where $\bar D_\omega$ is defined as \eqref{eq:Domega} with $\alpha_\o,\beta_\o$ replaced by $\bar\alpha_\o,\bar\beta_\o$.

Finally, we compute the contribution $b^{(0,2)}_{r,r'}(p)$ to the
Fourier transform of $\mathcal E_0(E_{x,r};E_{0,r'};V_4) $ where none
of the fields of $V_4$ are contracted among themselves.  First we write
in Fourier space
\begin{multline}
  V_4=-u\ik\ikp\ip \hat\psi^+_{k+p}\hat\psi^-_k\hat\psi^+_{k'-p}\hat\psi^-_{k'}\\
  \times(e^{i(k'_1+k'_2-p_2)}+e^{i(k'_1+k'_2-p_1)})\\
  =-\frac u4 \ik\ikp\ip \hat\psi^+_{k+p}\hat\psi^-_k\hat\psi^+_{k'-p}\hat\psi^-_{k'}W(k,k',p)
\end{multline}
where
\begin{multline}
  W(k,k',p)=
  e^{i(k_1'+k_2'-p_2)}+e^{i(k_1'+k_2'-p_1)} +e^{i(k_1+k_2+p_2)} +e^{i(k_1+k_2+p_1)}\\
  -e^{i(k_2+k'_1+p_2)}-e^{i(k_1+k'_2+p_1)}-e^{i(k_2'+k_1-p_2)}-e^{i(k_1'+k_2-p_1)}.\label{eqW}
\end{multline}
The second expression for $V_4$ is obtained by symmetrizing over the
four possible ways of ordering the fields $\hat\psi^+_{k+p}\hat\psi^-_k\hat\psi^+_{k'-p}\hat\psi^-_{k'}$ in such a way that the order of the upper indices is $(+,-,+,-)$.
Therefore, $b^{(0,2)}_{r,r'}(p)$  equals
\begin{eqnarray}
  -u K_r K_{r'}\ik\ikp W(k,k',p)\times\label{b02}\\
  \times\frac{e^{-i(k'-p)v_r-i(k+p)v_{r'}}}{\mu(k)\mu(k+p)\mu(k')\mu(k'-p)}.
\nonumber\end{eqnarray}
  Note that
\begin{eqnarray}
  \label{eq:W2}
  W(k,k',0)=2(e^{i k_1}-e^{i k'_1})(e^{i k_2}-e^{i k'_2}).
\end{eqnarray}
Then,  
\begin{eqnarray}
  \label{eq:dopbol}
  b^{(0,2)}_{r,r'}(p)= -2 u K_r K_{r'}\ik\ikp\times\\
  \times\frac{e^{-ik'v_r-ikv_{r'}}(e^{i k_1}-e^{i k'_1})(e^{i k_2}-e^{i k'_2})}{\mu(k)\mu(k+p)\mu(k')\mu(k'-p)}+R(p).
\end{eqnarray}
From this expression we want to extract the term that is discontinuous
at $p=0$.  Expanding the product in \eqref{eq:W2} we see that the
integral \eqref{eq:dopbol} is given by a sum of combinations of
integrals of the type $\mathcal I_{a}$ as in Lemma \ref{lemma:formulozzo}
below. Applying \eqref{eq:Ibis} to each of the integrals, it is easily
checked that the terms $(D_{-\o}(p)/D_\o(p))^2$ cancel and one is left
with\footnote{Here we assume, without loss of generality, that $\cos( p^+_1)>0$. Since $p^+_1+p^-_1=\pi$, we are just
  deciding which of the two zeros of $\mu(\cdot)$ we call $p^+$.  }
\begin{multline}
  \label{eq:dopbolfor}
  -\frac{i u K_r K_{r'}}{\pi(\alpha_+\beta_--\alpha_-\beta_+)}\sum_{\o=\pm}\frac{D_{-\o}(p)}{D_\o(p)}(e^{-i p^\o v_{r'}} U^\o_r+e^{-i p^\o v_{r}}U_{r'})+R(p),\\
U^\o_r=
e^{i p^\o(1,1)}C(-v_r)+C(-v_{r}+(1,1))\\-e^{i p^\o(1,0)}C(-v_r+(0,1))-e^{i p^\o(0,1)}C(-v_r+(1,0)).
\end{multline}
Using \eqref{eq:nota6} we see that
\begin{eqnarray}
  \label{eq:nota7}
  {(U^\o_r)^*}=U^{-\o}_r
\end{eqnarray}
and from  \eqref{eq:Ibis} one checks that
\begin{eqnarray}
  \nn
  U^+_1&=&\frac{e^{i p^+_2}}{2\pi}\int_{p^-_1}^{p^+_1+2\pi}d\theta \frac{e^{i p^+_1}-e^{i\theta}}{(K_1+K_2 e^{i\theta})^2}\\
  \label{eq:U1}
&&-\frac{2i}{\pi(\alpha_+\beta_--\alpha_-\beta_+)}\frac{\beta_+}{\beta_-}\cos(p^+_1)\cos(p^+_2)\\\nn
  U^+_2&=&\frac{e^{i p^+_2}}{2\pi}\int_{p^-_1}^{p^+_1+2\pi}d\theta \frac{e^{i\theta}(e^{i p^+_1}-e^{i\theta})}{(K_1+K_2 e^{i\theta})^2}\\
&&-\frac{2i}{\pi(\alpha_+\beta_--\alpha_-\beta_+)}e^{i p^-_1}\frac{\beta_+}{\beta_-}\cos(p^+_1)\cos(p^+_2)\\\nn
  U^+_3&=&\frac{1}{2\pi}\int_{p^+_1}^{p^-_1}d\theta \frac{e^{i\theta}(e^{i \theta}-e^{i p^+_1})}{(-K_3e^{i\theta}-K_4)^2}\\
&&-\frac{2i}{\pi(\alpha_+\beta_--\alpha_-\beta_+)}e^{i p^-_1+ip^-_2}\frac{\beta_+}{\beta_-}\cos(p^+_1)\cos(p^+_2)\\\nn
  U^+_4&=&\frac{1}{2\pi}\int_{p^+_1}^{p^-_1}d\theta \frac{e^{i \theta}-e^{i p^+_1}}{(K_3e^{i\theta}+K_4)^2}\\
&&-\frac{2i}{\pi(\alpha_+\beta_--\alpha_-\beta_+)}e^{i p^-_2}\frac{\beta_+}{\beta_-}\cos(p^+_1)\cos(p^+_2).
\end{eqnarray}

Summarizing, we obtained (dropping as usual the terms  $O(\l^2)$)
\begin{equation}
  \label{eq:cobinaz}
 \hat G^{(0,2)}_{r,r'}(p)=-\frac{i}{2\pi(\bar\alpha_+\bar\beta_--\bar\alpha_-\bar\beta_+)}\sum_{\o=\pm}\frac{\bar D_{-\o}(p)}{\bar D_\o(p)}\bar K_{\o,r }\bar K_{\o,r'}
e^{-i \bar p^\o(v_r+v_{r'})}+R(p)
\end{equation}
with
\begin{eqnarray}
  \label{eq:Ktilde}
  \bar K_{\o,r}:= K_r\left(e^{-i \bar p^\o v_r}-\frac{2\lambda t_r t_{r+2}}{K_{r}}W(p^\o)+2 u U^\o_r\right).
\end{eqnarray}
Thanks to \eqref{eq:symmab}, the prefactor of \eqref{eq:combinaz} is
real and \[\left(\frac{\bar D_-(p)}{\bar D_+(p)}\right)^*=\frac{\bar D_+(p)}{\bar D_-(p)}.\]
Finally, using \eqref{eq:nota0}, \eqref{eq:nota7} and \eqref{eq:symmp},
we see that the symmetry \eqref{eq:symmK} holds at first order in $\l$,
 so that  \eqref{eq:cobinaz} reduces to\begin{eqnarray}
  \label{eq:combinaz}
 \hat G^{(0,2)}_{r,r'}(p)=-\frac{i}{\pi(\bar\alpha_+\bar\beta_--\bar\alpha_-\bar\beta_+)}\Re\left[\frac{\bar D_{-}(p)}{\bar D_+(p)}\bar K_{+,r}
%e^{-i\bar p^+ v_r}
 \bar K_{+,r'}
    %     e^{-i \bar p^+v_{r'}}
  \right]+R(p).
\end{eqnarray}
As discussed in \cite[App. A]{GMT17b}, this asymptotic behavior of the
Fourier transform of the dimer-dimer correlation for $p\to0$ is
equivalent to the asymptotic behavior for large distances of the first
line of \eqref{dimdim} (the part proportional to $\bar H_{\o,r}$,
instead, in momentum space around $p=0$ can be absorbed into the error
term $R(p)$, because of the oscillating prefactor).
% It is easy to check that
% \begin{eqnarray}
%   \label{eq:4passi}
%   \sum_{r=1}^4 \bar K_r e^{-i \bar p^+ v_r}=O(\lambda^2).
% \end{eqnarray}
In order to prove at first order that the variance of the height difference grows logarithmically at large distances, as stated in Theorem \ref{th:2}, we 
need  to show that (see \cite[Theorem 2]{GMT19})
\begin{eqnarray}
  \label{eq:sespera}
  \frac{\Delta_2}{\bar \alpha_+}:=\frac{ \sum_{r\in\{1,4\}}\bar K_{+,r} % e^{-i \bar p^+ v_r}
  }{\bar \alpha_+}=
   \frac{\Delta_1}{\bar \beta_+}:= \frac{ \sum_{r\in\{1,2\}}\bar K_{+,r} % e^{-i \bar p^+ v_r}
  }{\bar \beta_+}
\end{eqnarray}
at order $\lambda$. This ratio is nothing but the first-order
expansion of $i \sqrt{A(\lambda)}$, cf. Theorem \ref{th:2}. Since for
$\lambda=0$ both ratios equal $i$ we write, at order $\l$,
\begin{align*}
  \bar\alpha_+&=\alpha_++2u \alpha_+^{(1)},\qquad\qquad\qquad \bar\beta_+=\beta_++2u \beta_+^{(1)},\\
  \Delta_1&=i \beta_++2u \Delta_1^{(1)},\qquad\qquad\quad\ \Delta_2=i\alpha_++2u\Delta^{(1)}_2,
\end{align*}
with $\alpha_+^{(1)}, \beta_+^{(1)},\Delta_1^{(1)},\Delta_1^{(2)}$
four $\lambda$-independent constants.  Then, we need to prove that
\begin{eqnarray}
  \label{eq:sp2}
  (\Delta^{(1)}_2-i\alpha_+^{(1)})\beta_+- (\Delta^{(1)}_1-i\beta_+^{(1)})\alpha_+=0.
\end{eqnarray}
From \eqref{a+b+} and \eqref{adm} it follows that
\begin{multline}
  \label{eq:a+b+bis}
  \alpha_+^{(1)}= -i g(v_1)e^{i(p^+_1+p^+_2)}+ig(v_4)e^{i p^+_1}+i \alpha_+ c^+_1 +
  c^+_2 t_3 e^{i (p^+_1+p^+_2)}\\
  \beta_+^{(1)}=-i g(v_1)e^{i(p^+_1+p^+_2)}+i g(v_2)e^{i p^+_2}+i \beta_+ c^+_2+c^+_1 t_3 e^{i (p^+_1+p^+_2)}
\end{multline} 
and from \eqref{eq:ptilde} and \eqref{eq:Ktilde} we see that
\begin{eqnarray}
  \label{eq:delta1}
  \Delta_2^{(1)}&=&-W(p^+)+K_1 U^+_1+K_4 U^+_4+i c^+_2 K_4 e^{i p^+_2}\\
  \Delta_1^{(1)}&=&-W(p^+)+K_1 U^+_1+K_2 U^+_2+i c^+_1 K_2 e^{i p^+_1}.
\end{eqnarray}
Therefore,
\begin{multline}
 \Delta_2^{(1)}-i \alpha_+^{(1)}=K_1 U^+_1+K_4 U^+_4-g(v_1)e^{i(p^+_1+p^+_2)}+g(v_4)e^{i p^+_1}\\
  \Delta_1^{(1)}-i \beta_+^{(1)}=K_1 U^+_1+K_2 U^+_2-g(v_1)e^{i(p^+_1+p^+_2)}+g(v_2)e^{i p^+_2}
\end{multline}
where we used $c_1^+\alpha_++c_2^+\beta_+=W(p^+)$ as in \eqref{eq:fpo}.
Then, using the definition of $U^+_r$, the l.h.s. of \eqref{eq:sp2} 
is
\begin{multline}
  \label{eq:cafebriosc}
  -\frac2\pi\frac{\beta_+}{\beta_-}\cos(p^+_1)\cos(p^+_2)\\-g(v_1)e^{i(p^+_1+p^+_2)}(\beta_+-\alpha_+)+\beta_+g(v_4)e^{i p^+_1}-\alpha_+g(v_2)e^{i p^+_2}\\
  -\alpha_+\frac{e^{ip^+_2}}{2\pi}\int_{p^-_1}^{p^+_1+2\pi}\frac{e^{i p^+_1}-e^{i\theta}}{K_1+K_2 e^{i\theta}}d\theta
  +\beta_+K_1\frac{e^{i p^+_2}}{2\pi}\int_{p^-_1}^{p^+_1+2\pi}\frac{e^{i p^+_1}-e^{i\theta}}{[K_1+K_2 e^{i\theta}]^2}d\theta\\
  +\beta_+ K_4\frac1{2\pi}\int_{p^+_1}^{p^-_1}\frac{e^{i\theta}-e^{i p^+_1}}{[K_3 e^{i\theta}+K_4]^2}d\theta.
\end{multline}
A simple application of the residue theorem shows that
\begin{eqnarray}
  \label{eq:Goo}
  g(v_1)=\frac1{2\pi}\int_{p^-_1}^{p^+_1+2\pi}\frac1{K_1+K_2 e^{i\theta}}d\theta\\
  g(v_2)=\frac1{2\pi}\int_{p^-_1}^{p^+_1+2\pi}\frac{e^{i\theta}}{K_1+K_2 e^{i\theta}}d\theta\\
  g(v_4)=\frac1{2\pi}\int_{p^+_1}^{p^-_1}\frac1{K_3e^{i\theta}+K_4}d\theta.
\end{eqnarray}
Therefore, the terms proportional to $\alpha_+$ in \eqref{eq:cafebriosc} cancel and one is left with
\begin{multline}
  \label{eq:Goo1}
   -\frac2\pi\frac{\beta_+}{\beta_-}\cos(p^+_1)\cos(p^+_2)\\-e^{i(p^+_1+p^+_2)}\beta_+\frac1{2\pi}\int_{p^-_1}^{p^+_1+2\pi}\frac1{K_1+K_2 e^{i\theta}}d\theta+\beta_+e^{i p^+_1}\frac1{2\pi}\int_{p^+_1}^{p^-_1}\frac1{K_3e^{i\theta}+K_4}d\theta\\
  +\beta_+K_1\frac{e^{i p^+_2}}{2\pi}\int_{p^-_1}^{p^+_1+2\pi}\frac{e^{i p^+_1}-e^{i\theta}}{[K_1+K_2 e^{i\theta}]^2}d\theta\\
  +\beta_+ K_4\frac1{2\pi}\int_{p^+_1}^{p^-_1}\frac{e^{i\theta}-e^{i p^+_1}}{[K_3 e^{i\theta}+K_4]^2}d\theta
  % =   -\frac2\pi\frac{\beta_+}{\beta_-}\cos(p^+_1)\cos(p^+_2)
  % +\beta_+\frac{K_4+K_3 e^{i p^+_1}}{2\pi}\int_{p^+_1}^{p^-_1}\frac{e^{i\theta}}
  % {[K_3 e^{i\theta}+K_4]^2}d\theta\\
  % -\beta_+(K_1+K_2 e^{i p^+_1})\frac{e^{i p^+_2}}{2\pi}\int_{p^-_1}^{p^+_1+2\pi}\frac{e^{i\theta}}{[K_1+K_2 e^{i\theta}]^2}d\theta
  \\
 % \label{cfbs}
  = -\frac2\pi\frac{\beta_+}{\beta_-}\cos(p^+_1)\cos(p^+_2) +\beta_+\frac{K_4+K_3 e^{i p^+_1}}{2\pi i }\int_{e^{i p^+_1}}^{e^{i p^-_1}}\frac{dz}
  {(K_3z+K_4)^2}\\
  -\beta_+(K_1+K_2 e^{i p^+_1})\frac{e^{i p^+_2}}{2\pi i}\int_{e^{i p^-_1}}^{e^{i p^+_1}}\frac{dz }{(K_1+K_2 z)^2}.
\end{multline}
Computing the integrals in the complex plane and using the explicit 
expressions  for  $\alpha_\o,\beta_\o$, which follow from their definition (see the statement of Theorem 1), one finds that
\eqref{eq:Goo1} is zero, so that \eqref{eq:sp2} holds.

In addition, one sees that the ratio \eqref{eq:sespera} is given by
\begin{equation}\label{sqrtA}
i \left(
1+\frac{4u}{\pi i} \frac{\cos(p^+_1)\cos(p^+_2)}{\alpha_+\beta_--\alpha_-\beta_+}
  \right)+O(\l^2).
  \end{equation}
Recalling that this ratio equals $i\sqrt{A(\l)}$ and that $\alpha_\o^*=-\alpha_{-\o},\beta^*_\o=-\beta_{-\o}$, we find that $A(\l)=1+a\l+O(\l^2)$, with $a$ as in \eqref{eqa}.

\subsection{First-order computation of $\nu(\l)$} By arguing like in
\cite[Appendix B.3]{GMT17b}, in order to compute the first order
contribution to $\nu$, it is enough to extract the most divergent
part, as $p\to p^{-}-p^+$, of the first order contribution to
$\hat G^{(0,2)}_{r,r'}(p)$. For definiteness, we assume, generically,
that $p^+-p^-\neq p^--p^+$ mod $2\pi$ (the complementary case, that corresponds to zero average tilt, can be
treated analogously). The coefficient $\nu_1$ in the expansion
$\nu(\l)=1+\nu_1\l+O(\l^2)$ can be read from:
%$$\hat G^{(0,2)}_{r,r'}(p^{-}-p^++q)= \frac{t_1^2}{4\pi^2}\frac1{\Delta}(-2\pi)\nu_1\l (\log|q|)^2+ \ldots,$$
\be \label{comp}\hat G^{(0,2)}_{r,r'}(p^{-}-p^++q)= \frac{-i\l\n_1{t_1^2}}{\p(\a_+\b_--\a_-\b_+)}(\log|q|)^2+ \ldots,\ee
where the dots indicate lower order terms in $\l$ or in $q$, as $q\to0$; i.e., they indicate terms that are either $O(\l^2)$, or less divergent than $(\log|q|)^2$ as $q\to 0$.  
 
By inspection (see the analogous discussion a few lines before \cite[eq.(B.29)]{GMT17b}), the only term that diverges like $(\log|q|)^2$ as $q\to 0$ comes from $b_{r,r'}^{(0,2)}$, see 
\eqref{b02}; setting, e.g., $r=r'=1$, 
\begin{eqnarray}&&\hat G^{(0,2)}_{1,1}(p^{-}-p^++q)=-u t_1^2 \ik \ikp W(k,k',p^--p^++q)\times\nonumber\\
&&\qquad \qquad \times \frac{1}{\mu(k)\mu(k+p^--p^++q)\mu(k')\mu(k'-p^-+p^+-q)}+ \ldots,\nonumber\end{eqnarray}
where the dots indicate lower order terms, as explained above. The dominant contribution to the right side comes 
from the region where $k$ is close to $p^+$ and $k'$ is close to $p^-$. By explicitly computing the dominant contribution to the integrand in this region we find: 
\begin{eqnarray}&&\hat G^{(0,2)}_{r,r'}(p^{-}-p^++q)=-\frac{u t_1^2}{16\pi^4} \int_{[-\e,\e]^2}\!\!\! ds \int_{[-\e,\e]^2} \!\!\!ds'\ W(p^+,p^-,p^--p^+) 
\times\nonumber\\
&&\qquad \qquad \times \frac{1}{D_+(s)D_-(s+q) D_-(s')D_+(s'-q)}+ \ldots,\nonumber\end{eqnarray}
where $\e$ is an arbitrary, small enough, positive number. 
Recalling \eqref{eqW} and the fact that $p^++p^-=(\pi,\pi)$, we find: 
$$W(p^+,p^-,p^--p^+) =-2(e^{ip_1^+}-e^{ip_1^-})(e^{ip_2^+}-e^{ip_2^-})=-8\cos(p_1^+)\,\cos(p_2^+).$$
Moreover, 
$$\int_{[-\e,\e]^2}\!\!\! ds\ \frac{1}{D_+(s)D_-(s+q)}=\frac{4i\pi}{(\a_+\b_--\a_-\b_+)}\log|q|+\ldots$$
Putting things together we find: 
$$\hat G^{(0,2)}_{r,r'}(p^{-}-p^++q)=-\frac{8u t_1^2\cos(p_1^+)\,\cos(p_2^+)}{\pi^2(\a_+\b_--\a_-\b_+)^2}(\log|q|)^2+\ldots$$
By comparing with \eqref{comp}, we find, as desired,
$$\nu_1 = \frac{8(u/\l)}{\pi i}\frac1{\a_+\b_--\a_-\b_+}\cos(p_1^+)\,\cos(p_2^+).$$

  \subsection{A useful integral formula }

% See details in contoresidui.pdf nella cartella dimeri/pendenza

 \begin{Lemma}\label{lemma:formulozzo}
   Let $a=(a_1,a_2)\in \mathbb Z^2$ and
\begin{eqnarray}
  \label{eq:I}
  \mathcal I_{a}(p)=\ik\frac{e^{i k a}}{\mu(k)\mu(k+p)}.
\end{eqnarray}
Assume that $\cos(p^+_1)>0$. Then, one has:
\begin{equation}
  \label{eq:Ibis}
 \mathcal  I_{a}(p)=\frac i{2\pi}\frac1{\alpha_+\beta_--\alpha_-\beta_+}\sum_{\omega=\pm}\frac{D_{-\omega}(p)}{D_\omega(p)}e^{ip^\omega a}
+C(a)+R_0(p)\end{equation}
where $R_0(p)$ tends continuously to zero as $p\to0$ and
\begin{gather}
C(a)=-\frac i{2\pi}\frac1{\alpha_+\beta_--\alpha_-\beta_+}\sum_{\omega=\pm}\frac{\beta_{-\omega}}{\beta_\omega}e^{ip^\omega a}\\
+\frac{(1-a_2)}{2\pi}\int_{0}^{2\pi}d\theta e^{i \theta a_1}\frac{(K_1+K_2 e^{i \theta})^{a_2-2}}{(-K_3e^{i \theta}-K_4)^{a_2}}\\\times
\left[1_{\{a_2\le 0\}} 1_{\{p_1^-\le \theta\le p_1^++2\pi\}}-1_{\{a_2\ge 1\}}1_{\{p_1^+\le \theta\le p_1^-\}}\right].
\end{gather}
 \end{Lemma}
The proof  is obtained with a straightforward although a bit lengthy application of the residue theorem; we skip details.
Note also that 
%\begin{eqnarray}
%  \label{eq:nota5}
 $ ({\mathcal I_{a}(p)})^*=(-1)^{a_1+a_2}\mathcal I_{a}(-p)$
%\end{eqnarray}
so that
\begin{eqnarray}
  \label{eq:nota6}  C(a)^*=(-1)^{a_1+a_2}C(a).
\end{eqnarray}

{\bf Acknowledgements.}  This review is based on a longstanding collaboration with Vieri Mastropietro, whom we thank for countless inspiring discussions. 
This work has been supported by the European Research Council (ERC) under the European Union's Horizon 2020 research and innovation programme 
(ERC CoG UniCoSM, grant agreement n.724939). F.T.  was  partially  supported  by
the  CNRS  PICS  grant  151933, by ANR-15-CE40-0020-03 Grant LSD, 
ANR-18-CE40-0033 Grant DIMERS and  by Labex MiLyon (ANR-10-LABX-0070). This work was started during a long-term stay of A.G. at Univ. Lyon-1, co-funded by Amidex and CNRS.

\end{document}